\documentclass[11pt,a4paper,twoside]{amsart}

\usepackage{amsmath,amssymb,amsfonts,amsthm,mathrsfs}
\usepackage{times,hyperref,color}
\usepackage[toc,page]{appendix}

\let\mathcal\mathscr

\newtheorem{The}{Theorem}[section]

\newtheorem{Theorem}{Theorem}[section]
\newtheorem{Claim}[The]{Claim}
\newtheorem{Proposition}[The]{Proposition}

\newtheorem{Lemma}[The]{Lemma}

\theoremstyle{definition}
\newtheorem{Definition}[The]{Definition}
\newtheorem{Convention}[The]{Convention}
\newtheorem{Remark}[The]{Remark}


\begin{document}


\title{
Totally nondegenerate models
\\
and standard manifolds
\\
in CR dimension one
}

\author{Masoud Sabzevari}
\address{Department of Mathematics,
University of Shahrekord, 88186-34141 Shahrekord, IRAN and School of
Mathematics, Institute for Research in Fundamental Sciences (IPM), P.
O. Box: 19395-5746, Tehran, IRAN}
\email{sabzevari@math.iut.ac.ir}

\date{\number\year-\number\month-\number\day}

\maketitle

\begin{abstract}
It is shown that two Levi-Tanaka and infinitesimal CR automorphism algebras, associated with a totally nondegenerate model of CR dimension one are isomorphic. As a result, the model surfaces are maximally homogeneous and standard. This gives an affirmative answer  in CR dimension one to a certain question formulated by Beloshapka.
\end{abstract}

\pagestyle{headings} \markright{Totally nondegenerate models and standard CR manifolds
}

\section{Preliminary definitions and results}
\label{sec-prel}

For an arbitrary smooth real manifold $M$, an even rank subbundle $T^cM\subset TM$ is called an {\sl almost CR structure} if it is equipped with a fiber preserving {\sl complex structure map} $J:T^cM\rightarrow T^cM$ satisfying $J\circ J=-id$. In this case, $M$ is called an {\sl almost CR manifold} of CR dimension $n:=\frac{1}{2}\cdot({\rm rank}\, T^cM)$ and codimension $k:={\rm dim}\, M-2n$. According to the principles in CR geometry (\cite{BER, 5-cubic}), the {\it complexified bundle} $\mathbb C\otimes T^cM$ decomposes as:
\[
\mathbb C\otimes T^cM:=T^{1,0}M\oplus T^{0,1}M,
\]
where:
\[
T^{1,0}M:=\Big\{X-i\,J(X): \ X\in T^cM\Big\}
\]
and $T^{0,1}M=\overline{T^{1,0}M}$. By definition, $M$ is called a {\sl CR manifold} with the {\sl CR structure} $T^cM$ if $T^{1,0}M$ is involutive in the sense of Frobenius. Such CR manifold is called a {\sl generic} submanifold of $\mathbb C^{n+k}$ if it can be represented locally as the graph of some $k$ defining functions $f_1,\ldots, f_k$ with $\partial f_1\wedge\ldots\wedge\partial f_k\neq 0$ ({\it cf.} \cite{BER}).

\subsection{Totally nondegenerate CR manifolds of CR dimension one}

Let $ M\subset\mathbb C^{1+k}$ be a real analytic generic submanifod of CR dimension one, codimension $k$, and hence of real dimension $2+k$. As is known (\cite{BER, Merker-Porten-2006, 5-cubic}), the holomorphic subbundle $T^{1,0}{ M}\subset\mathbb C\otimes T{ M}$ can be generated by some single holomorphic vector field $\mathcal L$. Set $D_1:=T^{1,0}{ M}+ T^{0,1}{ M}$ and also define successively $D_j=D_{j-1}+[D_1,D_{j-1}]$ for $j>1$. The iterated Lie brackets between the generators $\mathcal L$ and $\overline{\mathcal L}$ of $D_1$ induce a filtration:
 \[
D_1\subset D_2\subset D_3\subset\ldots
\]
 on the complexified tangent bundle $\mathbb C\otimes T{ M}$. Our distribution $D_1$ is {\sl minimal} or {\sl bracket-generating} if for each $p\in { M}$, there exists some (minimal) integer $\rho(p)$ satisfying $D_{\rho(p)}(p)=\mathbb C\otimes T_p{ M}$. Moreover, it is {\sl regular}, if the already mentioned function $\rho$ is constant. In this case, the number $\rho:=\rho(p)$ is called the {\sl degree of nonholonomy} of the distribution $D_1$.

\begin{Definition} (see \cite[Definition 1.1]{Bel-Conj}).
\label{def-totally-nondegenerate}
An arbitrary (local) real analytic CR generic submanifold ${ M}\subset\mathbb C^{1+k}$ of CR dimension one and codimension $k$ is {\sl totally nondegenerate} of the {\sl length} $\rho$ whenever the distribution $D_1=T^{1,0}{ M}+T^{0,1}{ M}$ is regular with the {\it minimum} possible degree of nonholonomy $\rho$. In this case, we have the induced filtration:
\begin{equation}
\label{filtration}
D_1\varsubsetneq D_2\varsubsetneq\ldots\varsubsetneq D_\rho=\mathbb C\otimes T{ M},
\end{equation}
of the minimum possible length.
\end{Definition}

The notion of total nondegeneracy has a close connection with the theory of complex {\it free Lie algebras} ({\it see} \cite{Merker-Porten-2006, Bel-Conj}). Clearly, in order to achieve the minimum length $\rho$ of the above filtration, each subbundle $D_{\ell}\setminus D_{\ell-1}$ must contain {\it maximum possible number} of independent iterated Lie brackets between $\mathcal L$ and $\overline{\mathcal L}$ of the length $\ell$, taking it into account that respecting the skew-symmetry and Jacobi identity is unavoidable. This number, that we denote it by ${\sf n}_\ell$, is computable by means of the well-known Witt's (inductive) formula ({\it cf.} \cite[Theorem 2.6]{Merker-Porten-2006}). By a careful inspection on the above definition, one verifies that the length $\rho$ of our $k$-codimensional totally nondegenerate CR manifold $M$ is in fact the smallest integer $\ell$ satisfying:
\begin{equation}
\label{rho}
{\rm rank}_{\mathbb R}(\mathbb C \otimes T{ M})=2+k \leq  {\sf n}_\ell.
\end{equation}
Moreover\,\,---\,\,as is the case with the free Lie algebras\,\,---\,\,{\it no linear relation exists between the iterated brackets of $\mathcal L$ and $\overline{\mathcal L}$ in the lengths $\leq\rho-1$, except those generated by skew-symmetry and Jacobi identity}. Notice that this rule entry into force until the length $\rho-1$ and, by contrast, one may encounter unpredictable treatments of iterated brackets in the lengths $\geqslant\rho$.

Thus, as a frame for the complexified bundle $\mathbb C\otimes TM$, the distribution $D_\rho$ is generated by the iterated brackets between $\mathcal L$ and its conjugation $\overline{\mathcal L}$ up to the length $\rho$. Following \cite{Bel-Conj}, let us show this frame by:
\begin{equation}
\label{frame}
\Big\{\mathcal L_{1,1}, \mathcal L_{1,2}, \mathcal L_{2,3}, \ldots, \mathcal L_{\rho,2+k}\Big\}
\end{equation}
where $\mathcal L_{1,1}:=\mathcal L$, $\mathcal L_{1,2}:=\overline{\mathcal L}$ and $\mathcal L_{\ell,i}$ is the $i$-th appearing independent vector field obtained as an iterated bracket of the length $\ell$.

 In \cite{Beloshapka2004}, Beloshapka showed that after appropriate weight assignment to the complex coordinates $(z, w_1, \ldots, w_k)$, every $k$-codimensional totally nondegenerate submanifold of $\mathbb C^{1+k}$ can be represented as the graph of some $k$ real analytic defining functions ({\it cf.} \cite[Theorem 1.1]{Bel-Conj}):
 \begin{equation}
 \label{general}
\aligned
\left\{
\begin{array}{l}
w_1=\Phi_1(z,\overline z)+{\rm O}([w_1]),
\\
\ \ \ \ \ \ \ \vdots
\\
w_k=\Phi_k(z,\overline z,\overline w)+{\rm O}([w_k]),
\end{array}
\right.
\endaligned
\end{equation}
where $[w_j]$ is the assigned weight to $w_j$ and where $\Phi_j$ is a weighted homogeneous complex-valued polynomial in terms of $z, \overline z, w_j$ and other complex variables $w_\bullet$ of the weights $[w_\bullet]<[w_j]$. Moreover, ${\rm O}(t)$ denotes some certain sum of monomials of the weights $> t$. In this case, the weight $[w_k]$ of the last variable $w_k$ is equal to the length $\rho$ of such CR manifold.
 Beloshapka also introduced:
  \begin{equation}
\label{model-def-eq}
\aligned
M:= \ \ \left\{
\begin{array}{l}
w_1=\Phi_1(z,\overline z),
\\
 \ \ \ \ \ \ \ \vdots
\\
w_k=\Phi_k(z,\overline z,\overline w),
\end{array}
\right.
\endaligned
\end{equation}
as a {\it model} of all $k$-codimensional totally nondegenerate manifolds, represented by \thetag{\ref{general}}. He also established a practical way to construct the associated defining polynomials $\Phi_\bullet$ ({\it see} \cite{Beloshapka2004} or \cite[$\S2$]{Bel-Conj}). These models are all homogeneous, of finite type and enjoy several other nice properties (\cite[Theorem 14]{Beloshapka2004}) that exhibit their significance.

 \begin{Convention}
 \label{convention}
  We stress that throughout this paper, we only deal with Beloshapka's totally nondegenerate CR models and, for the sake of brevity, we call them by {\sl "CR models"} or {\sl "models"}. From now on, we {\it fix the notation} $M$ for a certain totally nondegenerate model of CR dimension one, codimension $k$ and length $\rho$. Moreover, since we mainly utilize the results of \cite{Bel-Conj} and as is the case with that paper, we also assume that $\rho\geqslant 3$.
\end{Convention}

\subsection{Symbol algebra}

As we saw, the regular distribution $D_1=T^{1,0}M+T^{0,1}M$ of $M$ induces the filtration \thetag{\ref{filtration}} of the minimum length $\rho$. Set $\frak g_{-1}:=D_{1}$ and $\frak g_{-\ell}:=D_\ell\setminus D_{\ell-1}$ for $\ell>1$. By definition, $\frak g_{-\ell}$ is actually the vector space generated by all iterated Lie brackets between $\mathcal L$ and $\overline{\mathcal L}$ of the precise length $\ell$. As is shown in \cite[Proposition 3.2]{Bel-Conj} ({\it see also} Remark 3.3 of that paper) and independent of the choice of the points $p\in M$ that the above subdistributions $D_j$ are taken on them, the vector space:
\[
\frak g_-:=\frak g_{-\rho}\oplus\ldots\oplus\frak g_{-1}
\]
equipped with the standard Lie bracket of vector fields is essentially a unique $\rho$-th kind graded Lie algebra of dimension $2+k$, satisfying $[\frak g_{-i}, \frak g_{-j}]=\frak g_{-(i+j)}$. This algebra is {\sl fundamental}, that is: it can be generated by means of the Lie brackets between the elements of $\frak g_{-1}$. In this case, the regular distribution $D_1$ is called {\sl of constant type} $\frak g_-$ and, moreover, $\frak g_-$ is called by the {\sl symbol algebra} of $M$. We emphasize that the Lie algebra $\frak g_-$ and the distribution $D_\rho$ are actually two equal spaces of which the former is regarded as a Lie algebra while the latter is regarded as a frame for the complexified bundle $\mathbb C\otimes TM$.

\subsection{Lie algebras of infinitesimal CR automorphisms}

Let $\{\Gamma_{1,1}, \Gamma_{1,2}, \ldots, \Gamma_{\rho, 2+k}\}$ to be the (lifted) coframe, dual to \thetag{\ref{frame}}, of an arbitrary totally nondegenerate CR manifold of codimension $k$ which is biholomorphic (or CR-diffeomorphic) to $M$. Recently in \cite{Bel-Conj}, we have studied\,\,---\,\,by means of \'{E}lie Cartan's classical approach\,\,---\,\,the problem of biholomorphic equivalence to $M$ and found its associated {\it constant type structure equations} as:
 \begin{equation}
\label{structure-equations-after-prolongation}
\aligned
\left [
\begin{array}{l}
d\,\Gamma_{\ell,i}=(n_i\,\alpha+\tilde n_i\,\overline\alpha)\wedge\Gamma_{\ell,i}
+\sum_{\ell_1+\ell_2=\ell}\,{\sf c}^{i}_{j,n}\,\Gamma_{\ell_1,j}\wedge\Gamma_{\ell_2,n} \ \ \ \ {\scriptstyle (\ell\,=\,1\,,\,\ldots\,,\,\rho, \ \ i\,=\,1\,,\,\ldots\,,\,2+k)},
\\
d\alpha=0,
\\
d\overline\alpha=0,
\end{array}
\right.
\endaligned
\end{equation}
where ${\sf c}^i_{j,n}$s are some constant integers and where $n_i$ and $\tilde n_i$, visible among the expression of $d\Gamma_{\ell,i}$, are respectively the number of appearing $\mathcal L_{1,1}$ and $\mathcal L_{1,2}$ in constructing $\mathcal L_{\ell,i}$ as an iterated bracket of them. Moreover, $\alpha$ and $\overline\alpha$ are two certain Maurer-Cartan 1-forms added after prolongation steps of the method. Occasionally, it is possible to have these forms as real, {\it i.e.} $\alpha=\overline\alpha$, depending upon the CR model $M$, under study (\cite{Bel-Conj}).

As is known (\cite{5-cubic, Bel-Conj}), if the final structure equations of an equivalence problem to a certain $r$-dimensional smooth manifold $\bf M$ equipped with some lifted coframe $\{\gamma^1,\ldots,\gamma^r\}$ is of the constant type:

\[
d\gamma^k
=
\sum_{1\leqslant i<j\leqslant r}\,
c^k_{ij}\,\gamma^i\wedge\gamma^j
\ \ \ \ \ \ \ \ \ \ \ \ \
{\scriptstyle{(k\,=\,1\,\cdots\,r),}}
\]
then $\bf M$ is (locally) diffeomorphic to an $r$-dimensional Lie group $\sf G$, where its corresponding Lie algebra $\frak g$ has the basis elements $\{{\sf v}_1,\ldots, {\sf v}_r\}$\,\,---\,\,corresponding to $\{\gamma^1,\ldots,\gamma^r\}$\,\,---\,\,with the {\it structure constants}:
\[
\big[{\sf v}_i,{\sf v}_j\big]
=
-\sum_{k=1}^r\,c^k_{ij}\,{\sf v}_k
\ \ \ \ \ \ \ \ \ \ \ \ \
{\scriptstyle{(1\,\leqslant\,i\,<\,j\,\leqslant\,r)}}.
\]

We discovered in \cite{Bel-Conj} that the Lie algebra $\frak g$ associated with the final structure equations \thetag{\ref{structure-equations-after-prolongation}} is actually the desired Lie algebra $\frak{aut}_{CR}(M)$ of infinitesimal CR automorphisms of $M$. In order to realize the structure of this algebra through the above discussion, let us associate ${\sf v}_{\ell,i}$ to the 1-form $\Gamma_{\ell,i}$  for $\ell=1,\ldots,\rho$ and $i=1,\ldots,2+k$ and also associate ${\sf v_{0}}$ and ${\sf v}_{\overline 0}$ to $\alpha$ and $\overline\alpha$ as the basis elements of $\frak{aut}_{CR}(M)$. Clearly in the case that $\alpha=\overline\alpha$ we will have ${\sf v}_0={\sf v}_{\overline 0}$.

\begin{Proposition}
\label{prop-aut-graded}
(cf. \cite[Proposition 6.1]{Bel-Conj}). The Lie algebra $\frak{aut}_{CR}(M)$ is graded of the form:
\[
\frak{aut}_{CR}(M):=\underbrace{\frak g_{-\rho}\oplus\ldots\oplus\frak g_{-1}}_{\frak g_-}\oplus\, \frak g_0,
\]
satisfying $[\frak g_{-i},\frak g_{-j}]=\frak g_{-(i+j)}$, for $i,j=0,\ldots,\rho$, where $\frak g_-$ is (isomorphic to) the $(2+k)$-dimensional symbol algebra of $M$, where each homogeneous component $\frak g_{-\ell}$ is constructed by the basis elements ${\sf v}_{\ell,i}$ and where $\frak g_0$ is an Abelian Lie subalgebra of dimension either $1$ or $2$, generated by ${\sf v}_0$ and ${\sf v}_{\overline 0}$. The Lie brackets between these basis elements are determined by the constant type structure equations \thetag{\ref{structure-equations-after-prolongation}}.
\end{Proposition}

\begin{Definition}
 A graded Lie algebra $\frak g:=\bigoplus_{i\in\mathbb Z}\frak g_i$ is {\sl transitive} whenever $[{\sf x}_i, \frak g_-]\neq 0$ for each nonzero element ${\sf x}_i\in\frak g_i$ with $i\geqslant 0$. Also, it is {\sl nondegenerate} if $[{\sf x}_{-1}, \frak g_{-1}]\neq 0$ for each nonzero element ${\sf x}_{-1}\in\frak g_{-1}$.
\end{Definition}

\begin{Proposition}
\label{prop-transitive}
The Lie algebra $\frak{aut}_{CR}(M)=\frak g_-\oplus\frak g_0$ is nondegenerate and transitive.
\end{Proposition}

\proof
As is known ({\it cf.} \cite[p. 201]{Medori-Nacinovich-1997}), every prolongation of the symbol algebra associated with a Levi nondegenerate CR manifold is nondegenerate. Total nondegeneracy of $M$ implies its Levi nondegeneracy and hence $\frak{aut}_{CR}(M)$, as a prolongation of the symbol algebra $\frak g_-$, is nondegenerate.
 In order to prove that it is transitive, we have to check possible Lie brackets between the generators ${\sf v}_{1,1}, {\sf v}_{1,2}$ of $\frak g_{-1}$ and ${\sf v}_{0}, {\sf v}_{\overline 0}$ of $\frak g_0$ by looking for the wedge products between $\Gamma_{1,1}, \Gamma_{1,2}$ and $\alpha, \overline\alpha$ throughout the structure equations \thetag{\ref{structure-equations-after-prolongation}}. Such products exist only in the structure equations $d\Gamma_{1,1}$ and $d\Gamma_{1,2}$. First let us consider the case $\alpha\neq\overline\alpha$, where we have:
\[
\aligned
d\Gamma_{1,1}=\alpha\wedge\Gamma_{1,1} \ \ \ \ \ {\rm and} \ \ \ \ d\Gamma_{1,2}=\overline\alpha\wedge\Gamma_{1,2}.
\endaligned
\]
This implies that:
\begin{equation}
\label{bracket-g0-g1-1}
\aligned
& [{\sf v}_0, {\sf v}_{1,1}]=-{\sf v}_{1,1},  \ \ \ \ [{\sf v}_0, {\sf v}_{1,2}]=0,
\\
& [{\sf v}_{\overline 0}, {\sf v}_{1,1}]=0,  \ \ \ \ \ \ \ \ \ \ \ [{\sf v}_{\overline 0}, {\sf v}_{1,2}]=-{\sf v}_{1,2}.
\endaligned
\end{equation}
 For the case $\alpha=\overline\alpha$, the subalgebra $\frak g_0$ is generated by the single element ${\sf v}_0$ and the above equations give:
\begin{equation}
\label{bracket-g0-g1-2}
[{\sf v}_0, {\sf v}_{1,1}]=-{\sf v}_{1,1},  \ \ \ \ [{\sf v}_0, {\sf v}_{1,2}]=-{\sf v}_{1,2}.
\end{equation}
In any case, one observes that the Lie algebra $\frak{aut}_{CR}(M)$ is transitive, as was expected.
\endproof

\subsection{Tanaka prolongation and standard manifolds}

In \cite{Tanaka-main-1970}, Noboru Tanaka showed that associated with each finite dimensional fundamental graded algebra $\frak m:=\bigoplus_{-\mu\leqslant i\leqslant -1}\frak m_{i}$, there exists a unique, up to isomorphism, Lie algebra $\frak g(\frak m):=\bigoplus_{i\geqslant -\mu} \frak g^i(\frak m)$, satisfying:

\begin{itemize}
\item[(i)] $\frak g^i(\frak m)=\frak m_i$, for each $i=-\mu,\ldots,-1$.
\item[(ii)] $\frak g(\frak m)$ is transitive.
\item[(iii)] $\frak g(\frak m)$ is the maximal Lie algebra with the above two properties.
\end{itemize}
This algebra is known as the {\sl (full) Tanaka prolongation of} $\frak m$. He also established a practical method to construct successively the components $\frak g^i(\frak m)$ ({\it cf.} \cite{Alekseevsky-Spiro-2001,  Medori-Nacinovich-1997, Merker-Sabzevari-CEJM, Naruki-1970, Tanaka-main-1970}). In particular, the zero component $\frak g^0(\frak m)$ is the collection of all derivations ${\sf d}:\frak m\rightarrow\frak m$ that preserve the gradation, {\it i.e.} ${\sf d}(\frak m_{-1})\subset\frak m_{-1}$.

\begin{Definition}
A graded Lie algebra $\frak m:=\bigoplus_{i<0}\frak m_i$ is said to be {\sl pseudocomplex} (or {\sl CR}) if there exists some complex structure map $J:\frak m_{-1}\rightarrow \frak m_{-1}$ satisfying $J\circ J=-id$ and:
\[
\big[{\sf x}_{-1}, {\sf y}_{-1}\big]=\big[J({\sf x}_{-1}), J({\sf y}_{-1})\big], \ \ \ \ \ \ \ \ \textrm{for each} \ \ {\sf x}_{-1}, {\sf y}_{-1}\in\frak m_{-1}.
\]
\end{Definition}

 In the case that the graded fundamental algebra $\frak m:=\bigoplus_{-\mu\leqslant i\leqslant -1}\frak m_{i}$ is pseudocomplex, one defines as follows the so-called {\sl Levi-Tanaka} prolongation $\mathcal G(\frak m):=\bigoplus_{-\mu\geqslant i} \mathcal G^i(\frak m)$ of $\frak m$, essentially as a transitive subalgebra of $\frak g(\frak m)$: first, for each $i\leqslant -1$, set $\mathcal G^i(\frak m):=\frak m_i$. By definition, the zero component $\mathcal G^0(\frak m)$ is the collection of all derivations ${\sf d}\in\frak g^0(\frak m)$ that respect the associated complex structure map $J$, {\it i.e.}
 \begin{equation}
 \label{J-respect}
 {\sf d}(J({\sf x}_{-1}))=J({\sf d}({\sf x}_{-1})), \ \ \ \ \ \ \textrm{for each} \ \  {\sf x}_{-1}\in\frak m_{-1}.
  \end{equation}
The Lie bracket between two elements ${\sf d}\in\mathcal G^0(\frak m)$ and ${\sf x}\in\frak m_{i}$ is defined as $[{\sf d}, {\sf x}]:={\sf d}({\sf x})$. Assuming that the components $\mathcal G^{l'}(\frak m)$ are already constructed for any $l' \leqslant l - 1$, the $l$-th component $\mathcal G^l(\frak m)$ of the prolongation consists of $l$-shifted graded linear morphisms $\mathfrak m \to \mathfrak m \oplus \mathcal G^0(\frak m) \oplus \mathcal G^1(\frak m) \oplus \cdots
\oplus \mathcal G^{ l-1}(\frak m)$ that are derivations, namely:
\begin{equation*}
\mathcal G^l(\frak m)
=
\Big\{ {\sf d} \in
\bigoplus_{k\leqslant-1}\, {\sf Lin}\big(\mathcal G^k(\frak m),\,\mathcal G^{k+l}(\frak m)\big)
\colon
{\sf d}([{\sf y},\,{\sf z}]) =
[{\sf d}({\sf
y}),\,{\sf z}] + [{\sf y},\,{\sf d}({\sf z})],
\ \ \ \ \
\forall\, {\sf y},\,{\sf
z}\in\mathfrak m \Big\}.
\end{equation*}
Now, for ${\sf d} \in \mathcal G^k(\frak m)$ and ${\sf e} \in \mathcal G^l(\frak m)$, by induction on the integer $k + l \geqslant 0$, one defines
the bracket $[ {\sf d}, \, {\sf e} ] \in \mathcal G^{ k
+ l}(\frak m) \otimes \mathfrak m^*$ by:
\begin{equation*}
[{\sf d},\,{\sf e}]({\sf x}) = \big[[{\sf d},\,{\sf x}],\,{\sf e}\big]
+ \big[{\sf
d},\,[{\sf e},\,{\sf x}]\big] \ \ \ \ \ \ \ \ \ \ \text{\rm for}\ \ {\sf
x}\in\mathfrak m.
\end{equation*}

In the case that $\frak m:=\frak g_-$ is the symbol algebra of a certain manifold $\bf M$, then $\mathcal G(\frak m)$ is called the {\it Levi-Tanaka algebra of $\bf M$}.

For the Lie algebra $\frak{aut}_{CR}(M)=\frak g_-\oplus\frak g_0$, associated with our fixed CR model $M$, the Lie brackets between the basis elements ${\sf v}_{1,1}, {\sf v}_{1,2}$ of $\frak g_{-1}$ and ${\sf v}_0, {\sf v}_{\overline 0}$ of $\frak g_0$ are presented in \thetag{\ref{bracket-g0-g1-1}} and \thetag{\ref{bracket-g0-g1-2}} in two possible cases of $\alpha\neq\overline\alpha$ and $\alpha=\overline\alpha$. For some technical reasons, we substitute these basis elements with:
\[
\aligned
&{\rm for} \ \ \frak g_{-1}:\ {\sf x}:={\sf v}_{1,1}+{\sf v}_{1,2},  \ \ \ \ \ \ {\sf y}:=i\,({\sf v}_{1,1}-{\sf v}_{1,2}),
\\
&{\rm for} \ \ \frak g_0:
\left\{
\begin{array}{l}
{\sf d}:={\sf v}_{0}+{\sf v}_{\overline 0},   \ \ \ \  \ \ \ \ \ {\sf r}:=i\,({\sf v}_{0}-{\sf v}_{\overline 0}), \ \ \ \ \textrm{where} \ \ \alpha\neq\overline\alpha,
\\
{\sf d}:={\sf v}_{0}, \ \ \ \ \ \ \ \ \ \ \ \ \ \ \ \ \ \ \ \ \  \ \ \ \ \ \ \ \ \ \ \ \ \ \ \ \ \ \ \ \ \ \ \ \ \ \ \ \textrm{where} \ \ \alpha=\overline\alpha.
\end{array}
\right.
\endaligned
\]
Then, according to \thetag{\ref{bracket-g0-g1-1}} and \thetag{\ref{bracket-g0-g1-2}} we have:
\begin{equation*}
\label{bracket-g0-g1-main}
\left\{
\begin{array}{l}
\,[{\sf d}, {\sf x}]=-{\sf x}, \ \ \ [{\sf d}, {\sf y}]=-{\sf y}, \ \ \ [{\sf r}, {\sf x}]=-{\sf y}, \ \ \ [{\sf r}, {\sf y}]={\sf x} \ \ \  \textrm{where} \ \  \alpha\neq\overline\alpha,
\\
\\
\, [{\sf d}, {\sf x}]=-{\sf x}, \ \ \ [{\sf d}, {\sf y}]=-{\sf y}\ \ \  \ \ \ \ \ \ \ \ \ \ \ \ \ \ \ \ \ \ \ \ \ \ \ \ \ \ \ \ \ \ \ \ \ \ \ \ \ \ \ \ \ \ \textrm{where} \ \  \alpha=\overline\alpha.
\end{array}
\right.
\end{equation*}
Furthermore, let us define the complex structure map $J:\frak g_{-1}\rightarrow \frak g_{-1}$ by $J({\sf x})= {\sf y}$ and $J({\sf y})= -{\sf x}$.
\begin{Proposition}
\label{prop-aut-subset-levi-tanaka}
Assume as before that $\frak{aut}_{CR}(M)=\frak g_-\oplus\frak g_0$. By the already defined complex structure map $J$, the symbol algebra $\frak g_-$ of $M$ is pseudocomplex. Moreover, we have $\frak g_0\subseteq\mathcal G^0(\frak g_-)$.
\end{Proposition}

\proof
By the above definition of $J$, one readily verifies that for two basis elements $\sf x, y$ of $\frak g_{-1}$ we have $[{\sf x}, {\sf y}]=[J({\sf x}), J({\sf y})]$
which implies that $\frak g_-$ is pseudocomplex. For the second part of the assertion, first notice that according to Proposition \ref{prop-transitive}, $\frak{aut}_{CR}(M)$ is transitive and hence, by definition, it is a subalgebra of the Tanaka prolongation $\frak g(\frak g_-)=\bigoplus_{-\rho\geqslant i}\frak g^i(\frak g_-)$. Consequently we have $\frak g_0\subset\frak g^0(\frak g_-)$. Then it suffices to show that the elements of $\frak g_0$ respect the above complex structure map $J$. In other words, we have to show that $J([\mathcal D, \mathcal X])=[\mathcal D, J(\mathcal X)]$ for $\mathcal D={\sf d, r}$ and $\mathcal X= {\sf x, y}$ ({\it cf.}  \thetag{\ref{J-respect}}). It needs just some elementary computations that we leave them to the reader.
\endproof

\subsection{Main result} After providing preliminary definitions and results, concerning the subject, now we are ready to explain precisely the main aim of this paper. First we need the following two crucial definitions;

\begin{Definition}
({\it cf.} \cite{Medori-Nacinovich-1997}).
Let $\mathcal G(\frak m):=\bigoplus_{i\geqslant -\mu} \mathcal G^i(\frak m)$ to be the Levi-Tanaka prolongation of a pseudocomplex fundamental Lie algebra $\frak m$. Assume that $\sf G$ is the connected and simply connected Lie group with the Lie algebra $\mathcal G(\frak m)$ and also ${\sf G}_+$ is a closed analytic Lie subgroup of $\sf G$ with $\mathcal G_+(\frak m):=\bigoplus_{i\geqslant 0} \mathcal G^i(\frak m)$ as its Lie algebra. Then, the simply connected $\sf G$-homogeneous space $S(\mathcal G(\frak m)):=\frac{\sf G}{{\sf G}_+}$  is called the {\sl standard manifold} associated with the Levi-Tanaka prolongation $\mathcal G(\frak m)$.
\end{Definition}

\begin{Remark}
\label{rem}
 The above homogeneous manifold $S(\mathcal G(\frak m))=\frac{\sf G}{{\sf G}_+}$ is actually CR. To introduce its associated CR structure, consider first the natural projection $\pi:{\sf G}\rightarrow S(\mathcal G(\frak m))$ and let $\tt e$ and $o$ to be the identity elements of the groups $\sf G$ and $\frac{\sf G}{{\sf G}_+}$. The CR structure $T^cS(\mathcal G(\frak m))$ of $S(\mathcal G(\frak m))$ is defined  at the identity as:
 \[
 T^c_{o}\,S\big(\mathcal G(\frak m)\big):=\pi_{\ast}(\frak m_{-1}).
 \]
 Now, at each arbitrary point $\pi({\tt g})$ of $S(\mathcal G(\frak m))$, the desired CR structure is defined as the translation of the above fiber $T^c_{o}\,S\big(\mathcal G(\frak m)\big)$, from the identity $o$ to it, through the left  multiplication map $L_{\tt g}$, {\it i.e.}
 \[
  T^c_{\pi(g)}\,S\big(\mathcal G(\frak m)\big):=L_{g \ast}\big(T^c_{o}\,S(\mathcal G(\frak m))\big).
 \]
Roughly speaking, this CR structure is indeed the extension of $\frak m_{-1}$ to arbitrary points of  $S(\mathcal G(\frak m))$. It is invariant by the action of $\sf G$ on $S(\mathcal G(\frak m))$. For more details, we refer the reader to \cite[$\S4$]{Medori-Nacinovich-1997}.
\end{Remark}

An arbitrary CR manifold is standard if it is biholomorphic to a certain standard CR manifold. We have also the following\,\,---\,\,seemingly different but completely relevant\,\,---\,\,definition;

\begin{Definition}
({\it cf.} \cite{Medori-Nacinovich-2001}).
    An arbitrary Levi nondegenerate CR manifold $\bf M$ with the symbol algebra $\frak m$ is {\sl maximally homogeneous} if ${\rm dim}\, \frak{aut}_{CR}({\bf M})={\rm dim}\,\mathcal G(\frak m)$.
\end{Definition}

In \cite{Medori-Nacinovich-2001}, Medori and Nacinovich showed that {\it a nondegenerate CR manifold $\bf M$, regular of type $\frak m$, is maximally homogeneous if and only if it is biholomorphic to the associated standard manifold $S(\mathcal G(\frak m))$}.

 In 2004, Beloshapka formulated the question of whether his CR models are standard or not ({\it see} \cite[Question 2]{Beloshapka2004}). Our main aim in this paper is to answer this question affirmatively in CR dimension one. More precisely, the main result of this paper is as follows;

\begin{Theorem}
\label{theorem-main}
For each Beloshapka's totally nondegenerate model of CR dimension one, two associated Levi-Tanaka and infinitesimal CR automorphism algebras are isomorphic. As a result, such models are maximally homogeneous and standard.
\end{Theorem}

This result not only answers in part Beloshapka's  question but also provides infinitely many examples of standard manifolds with the certain known geometric-algebraic structures. We prove this theorem at the next section for CR models of the lengths $\rho\geqslant 3$ ({\it cf.} Convention \ref{convention}). By this upcoming proof, it remains only one model of which the correctness of the theorem should be proved in its case. It is nothing but the length two model $\mathbb H\subset\mathbb C^2$ of codimension one which is known as the {\it Heisenberg sphere} and is defined in coordinates $(z,w)$ of $\mathbb C^2$ as the graph of the following single polynomial equation:
\[
w-\overline w=2i\,z\overline z.
\]
Both the Lie algebra $\frak{aut}_{CR}(\mathbb H)$ and the Levi-Tanaka prolongation $\mathcal G(\frak g_-)$ associated to this exceptional model are computed explicitly in \cite[$\S\S 2, 3$]{Merker-Sabzevari-CEJM}, where we found them as two isomorphic $8$-dimensional graded algebras. Accordingly, the equality $\frak{aut}_{CR}(\mathbb H)= \mathcal G(\frak g_-)$ is plainly satisfied in this case.

\begin{Remark}
It is worth to emphasize that according to \cite[Proposition 3]{Beloshapka2004}, Beloshapka has proved (at least implicitly) that every totally nondegenerate CR model $M$ is diffeomorphic to the standard model associated with its symbol algebra, as two smooth manifolds. But, anyway, in order to conclude that $M$ is a standard manifold we need to have the already mentioned diffeomorphism to be CR. Unfortunately, Beloshapka's result does not bring directly such desired feature and our main objective in the next section is actually to prove it.
\end{Remark}

\section{Proof of Theorem \ref{theorem-main}}
\label{sec-proof}

For our length $\rho\geqslant 3$ CR model $M$,  Proposition \ref{prop-aut-subset-levi-tanaka} indicates that $\frak{aut}_{CR}(M)=\frak g_-\oplus\frak g_0$ is a subalgebra of the Levi-Tanaka algebra $\mathcal G(\frak g_-)$. In this section we prove the reverse inclusion $\frak{aut}_{CR}(M)\supseteq\mathcal G(\frak g_-)$.

In this case that the $(-1)$-component $\frak g_{-1}$ of the symbol algebra $\frak g_-=\bigoplus_{-\rho\leqslant i\leqslant -1}\frak g_i$ is of dimension two and according to last subsection 5.6 of \cite{Medori-Nacinovich-1997}, $\mathcal G^j(\frak g_0)$ is trivial for all $j\geqslant 1$. Thus, the Levi-Tanaka algebra $\mathcal G(\frak g_-)$ associated with our CR model $M$ is of the short form:
\[
\mathcal G(\frak g_-)=\frak g_-\oplus\mathcal G^0(\frak g_-).
\]
Consequently, our problem reduces to prove the inclusion $\frak{aut}_{CR}(M)\supseteq\mathcal G^0(\frak g_-)$.

Let ${\sf Aut}_{CR}(M)$ to be the connected and simply connected Lie group of all CR automorphisms of $M$, namely the collection of all automorphisms $h:M\rightarrow M$  satisfying $h_\ast(T^cM)=T^cM$. The associated Lie algebra to this finite dimensional group is $\frak{aut}_{CR}(M)$. We also denote by ${\sf Aut}_0(M)$ the connected isotropy subgroup of ${\sf Aut}_{CR}(M)$ at the origin. All automorphisms belonging to this subgroup are linear (\cite{Bel-Conj}) and its associated Lie algebra is $\frak g_0$. Finally, let ${\sf G_-}$ be the connected and simply connected Lie subgroup associated with the symbol algebra $\frak g_-$ of $M$. According to \cite[Proposition 3]{Beloshapka2004}, our CR model $M$ is an ${\sf Aut}_{CR}(M)$-homogeneous space and there exists a certain diffeomorphism:
\begin{equation}
\label{diff}
\Gamma: M \longrightarrow \frac{{\sf Aut}_{CR}(M)}{{\sf Aut}_0(M)}={\sf G_-} \ \ \ \ \ \ \ {\rm with} \ \  \ \ \Gamma(0)={\tt e},
\end{equation}
where $\tt e$ is actually the identity element of $\sf G_-$.

Let us denote by ${\sf Aut}_J(\frak g_-)$ the Lie group of all automorphisms of $\frak g_-$, preserving the gradation and respecting the complex structure map $J$ defined on the pseudocomplex algebra $\frak g_-$. According to Corollary 3, page 76 of \cite{Tanaka-main-1970} ({\it see also} \cite[Proposition 1.120]{Knap}), the associated Lie algebra  to this group is the zero component $\mathcal G^0(\frak g_-)$ of the Levi-Tanaka algebra of $M$. On the other, in this case that $\sf G_-$ is connected and simply connected, two automorphism Lie groups ${\sf Aut}({\sf G_-})$ and ${\sf Aut}(\frak g_-)$ are isomorphic through the map:
\[
\aligned
\Phi: {\sf Aut}({\sf G_-})&\longrightarrow {\sf Aut}(\frak g_-)
\\
f&\mapsto f_{\ast {\tt e}},
\endaligned
\]
where $f_{\ast \tt e}$ is the differentiation of $f:{\sf G_-}\rightarrow {\sf G_-}$ at the identity element $\tt e$  ({\it cf.} \cite{Hochschild-1952}). Let ${\sf Aut}_J({\sf G_-})\subset {\sf Aut}({\sf G_-})$ contains all automorphisms $f$ with $f_{\ast {\tt e}}\in{\sf Aut}_J(\frak g_-)$. Then, clearly we have:

\begin{Lemma}
Through the above isomorphism $\Phi$, two Lie groups ${\sf Aut}_J(\frak g_-)$ and ${\sf Aut}_J({\sf G_-})$ are isomorphic. As a result, the Lie algebra associated with ${\sf Aut}_J({\sf G_-})$ is $\mathcal G^0(\frak g_-)$.
\end{Lemma}

We aim to show that the Lie group ${\sf Aut}_J({\sf G_-})$ can be regarded as a subgroup of ${\sf Aut}_{CR}(M)$. As a result of this claim, we have:
\begin{equation}
\label{claim}
{\sf Lie}\big({\sf Aut}_J({\sf G_-})\big)\subseteq {\sf Lie}\big({\sf Aut}_{CR}(M)\big) \ \ \ \textrm{or equivalently} \ \ \
\mathcal G^0(\frak g_-)\subseteq\frak{aut}_{CR}(M),
\end{equation}
 as was desired. For this purpose and denoting by ${\sf Aut}(M)$ the collection of all (not necessarily CR) automorphisms from $M$ to itself, we define:
 \begin{equation}
\label{Psi}
 \aligned
 \Psi:{\sf Aut}_J({\sf G_-})&\longrightarrow {\sf Aut}(M)
   \\
   f&\mapsto\Gamma^{-1}\circ f \circ \Gamma=:F
\endaligned
    \end{equation}
where $\Gamma$ is the above diffeomorphism \thetag{\ref{diff}}. Having $f$, as an automorphism of $\sf G_-$, and $\Gamma$ as two certain diffeomorphisms then $\Psi(f)$, that we denote it by $F$ henceforth, is an automorphism.
  We claim that;

 \begin{Claim}
 \label{claim-main}
 $F$ is a CR automorphism and hence belongs to ${\sf Aut}_{CR}(M)$.
 \end{Claim}

In order to prove this claim, we need first the following auxiliary lemma;

 \begin{Lemma}
 \label{lem-left-invariant}
There exists two independent basis vector fields ${\sf L}_1$ and ${\sf L}_2$ of $T^cM$ such that both ${\sf X}_1:=\Gamma_\ast({\sf L}_1)$ and ${\sf X}_2:=\Gamma_\ast({\sf L}_2)$, as two vector fields defined on $\sf G_-$, are left invariant.
 \end{Lemma}

 \proof
 Let ${\sf L}_1$ and ${\sf L}_2$ be two arbitrary generators of $T^cM$ and consider ${\sf X}_1:=\Gamma_\ast({\sf L}_1)$ and ${\sf X}_2:=\Gamma_\ast({\sf L}_2)$ as two, {\it not necessarily left invariant}, vector fields defined on $\sf G_-$. As is known (\cite[Proposition 4(c)]{Beloshapka2004}), the Lie subalgebra $\frak g_{-1}$ of $\frak{aut}_{CR}(M)$ can be generated by ${\sf x}_1:={\sf X}_1({\tt e})$ and ${\sf x}_2:={\sf X}_2({\tt e})$. Now, assume that $\widetilde{\sf X}_1$ and $\widetilde {\sf X}_2$ are the (unique) independent left invariant vector fields associated with ${\sf x}_1$ and ${\sf x}_2$, respectively. Let ${\frak D}_1$ to be the distribution on $\sf G_-$ generated by $\widetilde{\sf X}_1$ and $\widetilde {\sf X}_2$ and define successively ${\frak D}_j:={\frak D}_{j-1}+[{\frak D}_1, {\frak D}_{j-1}]$ for $j>1$. Since $\widetilde{\sf X}_1$ and $\widetilde {\sf X}_2$ are left invariant, the value of each iterated bracket between them at an arbitrary point ${\tt g}\in \sf G_-$ is actually the translation of this value at $\tt e$ through the differential map $L_{{\tt g} \ast}$ of the left multiplication by $\tt g$, {\it i.e.}
 \[
 \big[\widetilde{\sf X}_{i_1}, [\widetilde{\sf X}_{i_2}, [\widetilde{\sf X}_{i_3}, [\ldots, [\widetilde{\sf X}_{i_{\ell-1}}, \widetilde{\sf X}_{i_\ell}]]]]\big]_{\tt g}=L_{{\tt g} \ast}\Big( \big[{\sf x}_{i_1}, [{\sf x}_{i_2}, [{\sf x}_{i_3}, [\ldots, [{\sf x}_{i_{\ell-1}}, {\sf x}_{i_\ell}]]]]\big]\Big), \ \ \ \ \ \ {\scriptstyle (i_j\,=\, 1\,,\,2)}.
 \]
 This indicates that at each arbitrary point of $\sf G_-$, the filtration ${\frak D}_1\subset{\frak D}_2\subset {\frak D}_3\subset\ldots$
has exactly similar treatment as that at the identity element $\tt e$. Consequently, at each arbitrary point ${\tt g}\in\sf G_-$, the above filtration induces the same fundamental symbol algebra $\frak g_-$ generated by ${\sf x}_1$ and ${\sf x}_2$.  Set $\widetilde{\sf L}_i:=\Gamma^{-1}_\ast( \widetilde{\sf X}_i)$ for $i=1,2$. The subbundle $\widetilde{T}^cM$ generated by these two independent vector fields can be regarded as a certain CR structure for $M$ via the complex structure map $\widetilde{J}:\widetilde{T}^cM\rightarrow \widetilde{T}^cM$ defined by $\widetilde{J}(\widetilde{\sf L}_1):=\widetilde{\sf L}_2$ and $\widetilde{J}(\widetilde{\sf L}_2):=-\widetilde{\sf L}_1$. According to Theorem 3, page 109 of \cite{BER}, determining a basis for the CR structure of an smooth generic submanifold depends only upon its defining functions and hence $\widetilde{T}^cM$ and $T^cM$ are two equivalent CR structures for $M$. Therefore, we can consider $\widetilde{\sf L}_1$ and $\widetilde{\sf L}_2$ as the basis elements of the CR structure of $M$ where their image vector fields $\widetilde{\sf X}_1=\Gamma_\ast(\widetilde{\sf L}_1)$ and $\widetilde{\sf X}_2=\Gamma_\ast(\widetilde{\sf L}_2)$ are left invariant on $\sf G_-$.
\endproof

\noindent
{\it Proof of Claim \ref{claim-main}.} Since $F=\Psi(f)$ is a diffeomorphism, then it suffices to show that $F_\ast(T^cM)=T^cM$. Let ${\sf L}_1$ and ${\sf L}_2$ be two basis vector fields of $T^cM$ such that ${\sf X}_1=\Gamma_\ast({\sf L}_1)$ and ${\sf X}_2=\Gamma_\ast({\sf L}_2)$ are left invariant, with ${\sf X}_1(e)={\sf x}_1$ and ${\sf X}_2(e)={\sf x}_2$ as the basis elements of $\frak g_{-1}$. By definition, since $f\in{\sf Aut}_J({\sf G_-})$ then its differentiation $f_{\ast e}$ preserves the gradation and hence we have:
\[
f_{\ast e}({\sf x}_1):=a_1\,{\sf x}_1+ a_2\,{\sf x}_2, \ \ \ {\rm and} \ \ \ f_{\ast e}({\sf x}_2):=b_1\,{\sf x}_1+ b_2\,{\sf x}_2,
\]
for some constant integers $a_\bullet$ and $b_\bullet$. Now, since ${\sf X}_1$ and ${\sf X}_2$ are left invariant, then it yields that:
\[
f_{\ast}({\sf X}_1):=a_1\,{\sf X}_1+ a_2\,{\sf X}_2, \ \ \ {\rm and} \ \ \ f_{\ast }({\sf X}_2):=b_1\,{\sf X}_1+ b_2\,{\sf X}_2,
\]
at each arbitrary point ${\tt g}\in\sf G_-$. Consequently, $f_\ast$ preserves as well the subdistribution  $\frak D_1:=\langle{\sf X}_1, {\sf X}_2\rangle$, {\it i.e.}
\[
f_\ast(\frak D_1({\tt g}))=\frak D_1({\tt h}), \ \ \ \ {\rm with} \ \ \ {\tt h}:=f({\tt g}).
\]
Let $p\in M$, with $\Gamma(p)={\tt g}$ and $f({\tt g})={\tt h}$. Since $\Gamma_\ast(T^cM)=\frak D_1$, then we have:
\[
\aligned
F_\ast(T^c_pM)&=\Gamma_{\ast {\tt h}}^{-1}\circ f_{\ast {\tt g}}\circ \Gamma_{\ast p}(T^c_pM)
\\
&=\Gamma_{\ast {\tt h}}^{-1}\circ f_{\ast {\tt g}}(\frak D_1({\tt g}))=\Gamma_{\ast {\tt h}}^{-1}(\frak D_1({\tt h}))
\\
&=T^c_{F(p)}M.
\endaligned
\]
This completes the proof. \ \ \ \ \ \ \ \ \ \ \ \ \ \ \ \ \ \ \ \ \ \ \ \ \ \ \ \ \ \ \ \ \ \ \ \ \ \ \ \ \ \ \ \ \ \ \  \ \ \ \ \ \ \ \ \ \ \ \ \ \ \ \ \ \ \ \ \ \ \ \ \ \ \ \ \ \ \ \ \ \ \ \ \ \ \  \ \ \ \ \ \ \ \ \ \ \ \ \ \ \ \ \ \ \ \ \ \ \ \ \ \ \ \ \ \ \ \ \ \ \ \ \ \ \  \ \ \ \ \ \ \ \  \ \ \ \ \ \ \ \ \ \ \ \ \ \ \ \ \ \ \ \ \ \ \ \ \ $\Box$

Consequently, we can take the map $\Psi$, introduced in \thetag{\ref{Psi}} as:
\[
\Psi:{\sf Aut}_J({\sf G_-})\longrightarrow {\sf Aut}_{CR}(M).
\]
One also readily verifies that this map is actually an injective Lie group homomorphism. Therefore, we can regard ${\sf Aut}_J({\sf G_-})$  as a Lie subgroup of ${\sf Aut}_{CR}(M)$. Then according to \thetag{\ref{claim}}, we have $\mathcal G^0(\frak g_-)\subseteq\frak{aut}_{CR}(M)$, as was desired. This completes the proof of the main Theorem \ref{theorem-main}.  \ \ \ \ \ \ \ \ \ \ \ \ \ \ \ \ \ \ \ \ \ \ \ \ \ \ \ \ \ \ \ \ \ \ \ \ \ \ \ \ \ \ \ \ \ \ \ \ \ \ \ \ \ \ \ \ \ \  \ \ \ \ \ \ \ \ $\Box$

\begin{Remark}
It may be worth to notice that since all the automorphisms $f\in {\sf Aut}_J({\sf G_-})$ fix the identity element $e$ of $\sf G_-$, then, its associated CR automorphism $F\in {\sf Aut}_{CR}(M)$ preserves the origin:
\[
F(0)=\underline{\Gamma^{-1}\circ \underline{f \circ \underline{\Gamma (0)}_{\ {\tt e}}}_{\ {\tt e}}}_{\ 0}=0.
\]
Therefore, the induced function $F$ belongs to the isotropy subgroup ${\sf Aut}_0(M)$ of ${\sf Aut}_{CR}(M)$ at the origin. Then, in a more precise manner, we can regard ${\sf Aut}_J({\sf G_-})$ as a subgroup of the isotropy group ${\sf Aut}_0(M)$. This implies the inclusion $\mathcal G^0(\frak g_-)\subseteq\frak g_0$ and consequently we have $\mathcal G^0(\frak g_-)=\frak g_0$ as was expected.
\end{Remark}

\subsection*{Acknowledgment}
The author expresses his sincere thanks to Mauro Nacinovich and Andrea Spiro for their
helpful comments during the preparation of this paper. Also he would like to thank Valerii Beloshapka for his helpful discussions about the novelty of the problem, considered in this paper.
The research of the author was supported in part by a grant from IPM, No. 94510061.

\end{document}